\documentclass[a4paper]{article}
\usepackage[english]{babel}
\usepackage[utf8]{inputenc}
\RequirePackage{amsmath,amsthm,amssymb,enumitem,graphicx,ifpdf}
\RequirePackage[left=1in,right=1in,top=1in,bottom=1in]{geometry}
\newtheorem{theorem}{Theorem}[subsection]

\newtheorem*{claim*}{Claim}
\newtheorem*{theorem*}{Theorem}

\newtheorem*{corollary*}{Corollary}

\newtheorem*{lemma*}{Lemma}
\theoremstyle{definition}

\newtheorem*{definition*}{Definition}

\newtheorem*{proposition*}{Proposition}

\newtheorem*{notation*}{Notation}

\newtheorem*{remark*}{Remark}

\numberwithin{equation}{section}
\newcommand{\R}{\mathbb{R}}

\newcommand{\N}{\mathbb{N}}

\usepackage{cite}

\title{A Simple Multiple Integral Solution to the Broken Stick Problem}
\author{Vivek Kaushik}
\begin{document}
\maketitle
\begin{abstract}
Fix an integer $n \ge 2.$ Partition the interval $[0,1]$ into $n+1$ different intervals $I_1, \ \dots \ , I_{n+1}.$ We solve the Broken Stick problem, which is to find the probability that the lengths of these $n+1$ intervals are valid side lengths of a polygon with $n+1$ sides. We will show that this probability is equal to $1-\frac{n+1}{2^{n}}$ using multiple integration.
\end{abstract}

\section{Introduction}
Partition the interval $[0,1]$ into $n+1$ different intervals $I_1, \ \dots \ , I_{n+1},$ where $n \ge 2$ is an integer. The \textit{Broken Stick problem} asks to find the probability these $n+1$ intervals have lengths that are valid side lengths of a polygon with $n+1$ sides. By classical geometry, this is equivalent to finding the probability that for each $1 \le i \le n+1,$ the length of $I_i$ is less than the sum of the lengths of $I_j$  for all $j \ne i.$

The answer to the Broken Stick problem is:
\begin{align}\label{Broken Stick Solution}
1-\frac{n+1}{2^{n}}.
\end{align}
The particular case $n=2$ is well-known. In this case, the Broken Stick problem asks to find the probability that the lengths of $I_1, I_2,$ and $I_3$ are valid side lengths of a triangle. The answer is $1/4,$ and multiple solutions from different areas of mathematics have appeared in the literature (see all of \cite{Goodman,Lemoine1,Lemoine2,Seneta,Poincare}). D'Andrea and Gomez \cite{Gomez} published a solution to the general case using combinatorial topology.

In this paper, we will give an alternative solution to the Broken Stick problem using multiple integration. We let $X_1, \ \dots \ , X_n$ be $n$ independent and identically distributed uniform random variables on $[0,1].$ Regarding $[0,1]$ as a stick, these random variables represent the $n$ breaking points of the stick. We make the ordering assumption that $X_1 \le \ \dots \ \le X_n$ and show that the solution to the Broken Stick problem with this ordering assumption is
\begin{align}\label{Broken Stick Solution with Ordering}
\frac{2^n-(n+1)}{2^{n} \ n!}.
\end{align} 
To do this, we formulate the necessary conditions on these random variables for them to satisfy the Broken Stick problem. Then, for each $1 \le k<n,$ we further impose the constraint that exactly $k$ of these random variables must be less than $1/2,$ and we solve the Broken Stick problem with this additional constraint by explicitly evaluating a multiple integral. Using a well-known combinatorial formula, we sum the results over all possible aforementioned $k$ and obtain the solution to the Broken Stick problem with only our original ordering assumption, namely \eqref{Broken Stick Solution with Ordering}. Finally, since there are $n!$ total orderings of $X_1, \ \dots \ ,X_n,$ multiplying the result in \eqref{Broken Stick Solution with Ordering} by $n!$ recovers the solution to the original Broken Stick problem without any ordering assumptions, namely \eqref{Broken Stick Solution}.

\section{Preliminaries}
In this section, we prove two formulas that are crucial for solving the Broken Stick problem. 

The first formula is an integration identity of a certain monomial function over an $n$-dimensional simplex.
\begin{theorem}[Integration Identity of Monomial Function Over a Simplex]\label{Integration Identity of Monomial Function Over A Simplex}

Suppose $p \in \N \cup \lbrace 0 \rbrace$ and suppose $a,b \in \R$ with $a \le b.$ For each $n \in \N,$ we have
\begin{align}\label{Integral of Monomial Function Over A Simplex}
\int_{a}^b \int_{0}^{x_1} \ \dots \ \int_{0}^{x_{n-1}}  x_n^p \ dx_n \ \dots \ dx_2 \ dx_1 &= \frac{b^{p+n} -a^{p+n}}{(p+1) \ \dots \ (p+n)}.
\end{align}
\end{theorem}

\begin{proof}
We proceed by induction on $n.$

In the case $n=1,$ the left hand side of \eqref{Integral of Monomial Function Over A Simplex} becomes the integral
\begin{align}\label{Case n=1 Integral}
\int_{a}^b x_1^p \ dx_1.
\end{align}
By the elementary power rule of integration, the integral in \eqref{Case n=1 Integral} is equal to $$\frac{b^{p+1} - a^{p+1}}{p+1},$$ which clearly agrees with the right hand side of \eqref{Integral of Monomial Function Over A Simplex} upon substituting $n=1.$

Now, suppose the claim holds for the case $n=k$ for some $k>1.$ We show that it must also hold for the case $n=k+1.$  In the case $n=k+1,$ the left hand side of \eqref{Integral of Monomial Function Over A Simplex} is equal to
\begin{align}\label{Case n=k+1 Integral}
\int_{a}^b \left( \int_{0}^{x_1} \ \dots \ \int_{0}^{x_{k}}  x_{k+1}^p \ dx_{k+1} \ \dots \ dx_2 \right) \ dx_1 &= \int_{a}^b \frac{x_1^{p+k}}{(p+1) \ \dots \ (p+k)} \ dx_1 \\
&= \frac{b^{p+k+1} - a^{p+k+1}}{(p+1) \ \dots \ (p+k+1)} \nonumber,
\end{align}
where we obtained the right hand side of \eqref{Case n=k+1 Integral} by using the induction hypothesis to evaluate the $k$-dimensional integral within the parentheses on the left hand side of \eqref{Case n=k+1 Integral}.
\end{proof}

The next formula is a well-known combinatorial formula.
\begin{theorem}[Combinatorial Identity]\label{Combinatorial Identity}
For all $n \in \N,$ we have 
\begin{align}\label{2^n identity}
\sum_{p=0}^n \binom{n}p &=2^n,
\end{align}
where 
\begin{align}\label{Binomial Coefficient}
\binom{n}p=\frac{n!}{(n-p)! \ p!}
\end{align}
is the binomial coefficient. 
\end{theorem}

\begin{proof}
The left hand side of \eqref{2^n identity} is precisely the expansion of $(1+1)^n$ using the Binomial Theorem. On the other hand, we have $(1+1)^n=2^n,$ which is the right hand side of \eqref{2^n identity}.
\end{proof}

\section{Multiple Integration Solution to the Broken Stick Problem}
We now are ready to solve the Broken Stick problem. As stated in the introduction, we let $n \ge 2$ be an integer and let $X_1, \ \dots \ , X_n$ be $n$ independent and identically distributed uniform random variables on $[0,1].$ Regarding $[0,1]$ as a stick, these random variables represent the breaking points of $[0,1]$ so that we obtain $n+1$ different intervals $I_1, \ \dots \ , I_{n+1}$ that partition $[0,1].$    The Broken Stick problem asks us to find the probability that the length of each $I_i$ is less than the sum of the lengths of $I_j$ for all $j \ne i.$ 

We first explicitly rewrite this probability in terms of conditions on $X_1, \ \dots \ , X_n$ assuming an ordering of these random variables. 
\begin{theorem}[Broken Stick Problem Assuming Ordering]\label{Broken Stick Problem Assuming Ordering}
Assume that $X_1 \le \ \dots \ \le X_n.$ Then, the Broken Stick problem is equivalent to computing the probability that these random variables satisfy the following conditions:
\begin{align} \label{Broken Stick Conditions Assuming Ordering}
X_1<\frac{1}{2}, \ X_1 \le X_{2} < \frac{1}{2} + X_1, \ \dots \ , \  X_{n-1} \le X_{n} < \frac{1}{2} + X_{n-1}, \ X_n > \frac{1}{2}.
\end{align}
\end{theorem}

\begin{proof}
Since $X_1 \le \ \dots \ \le X_n,$ the interval $[0,1]$ is partitioned into $I_1, \ \dots \ , I_{n+1},$ where
\begin{align*}
I_1 &= [0,X_1] \\
I_i &= [X_{i-1}, X_i], \quad 1<i \le n \\
I_{n+1} &= [X_{n}, 1].
\end{align*}
We let $\ell(I)$ denote the length of the interval $I.$ 

The condition that the length of $I_1$ is less than the sum of the lengths of all $I_j$ for $j \ne 1$ is equivalent to the condition
\begin{align*}
X_1 = \ell(I_1) &< \sum_{j \ne 1}^{n+1} \ell(I_j) \\
&=\ell(I_{n+1}) + \sum_{j \ne 1}^n \ell(I_j) \\
&= 1-X_n + \sum_{j=2}^n (X_{j}-X_{j-1}) \\
&= 1-X_1.
\end{align*}
Rearranging yields $X_1<1/2,$ which is precisely the first condition listed in \eqref{Broken Stick Conditions Assuming Ordering}. 

Let $1<i \le n.$ The condition that the length of $I_i$ is less than the sum of the lengths of all $I_j$ for $j \ne i$ is equivalent to the condition
\begin{align*}
0 \le X_{i}-X_{i-1} = \ell(I_i) &< \sum_{j \ne i}^{n+1} \ell(I_j) \\
&= \ell(I_1) + \ell(I_{n+1}) + \sum_{j \ne 1,n+1,i}^{n} \ell(I_j) \\
&= X_1 + 1-X_n + \sum_{j \ne 1,n+1,i}^{n} (X_j-X_{j-1}) \\
&= 1-X_i+X_{i-1}.
\end{align*}
Rearranging yields $X_i<1/2+X_{i-1}.$ On the other hand, by the ordering assumption, we have $X_{i-1} \le X_{i}.$ These two inequalities are the same as the $i$-th condition listed in \eqref{Broken Stick Conditions Assuming Ordering}.

Lastly, the condition that the length of $I_{n+1}$ is less than the sum of the lengths of all $I_j$ for $j \ne {n+1}$ is equivalent to the condition
\begin{align*}
1-X_{n} = \ell(I_{n+1}) &< \sum_{j \ne n+1}^{n+1} \ell(I_j) \\
&= \sum_{j \ne 1}^n \ell(I_j) \\
&= \sum_{j=2}^n (X_{j}-X_{j-1}) \\
&= X_n.
\end{align*}
Rearranging yields $X_n>1/2,$ which is precisely the last condition listed in \eqref{Broken Stick Conditions Assuming Ordering}.
\end{proof}

In addition to the ordering assumption we made on the random variables in \textbf{Theorem \ref{Broken Stick Problem Assuming Ordering}}, we further impose a condition that a certain number of these random variables must be less than or equal to $1/2.$ This will allow us to solve the Broken Stick problem by explicitly evaluating a multiple integral.

\begin{theorem}[Broken Stick Solution With Ordering and Extra Bound Assumptions]\label{Broken Stick Solution With Ordering and Extra Bound Assumptions}
As before, assume $X_1 \le \dots \ \le X_n.$ Fix $1 \le k <n,$ and let $p_{n,k}$ denote the probability that $X_1, \ \dots \ , X_n$ satisfy all the conditions listed in \eqref{Broken Stick Conditions Assuming Ordering} and the additional conditions $X_1, \ \dots \ , X_k<1/2, \ X_{k+1} \ge 1/2.$ 
Then, 
\begin{align}\label{p_nk Formula}
p_{n,k} &=\frac{\binom{n}k  -1}{2^n \ n!}.
\end{align}
\end{theorem}

\begin{proof}
To prove this, we set up and evaluate a multiple integral corresponding to $p_{n,k}.$

The integrand for this multiple integral will be the joint density function of $X_1, \ \dots \ , X_n,$ which is $1.$ 

We wish to set up the multiple integral so that we integrate in the order of $x_n, \ \dots \ , x_1.$ We determine the integral bounds corresponding to $x_n, \ \dots \ , x_1$ in that order. We obtain each integral bound for $x_i$ by obtaining a bound for the corresponding random variable $X_i.$ 

For each $k+1 <i \le n,$ we have $X_i \ge X_{i-1}$ and $X_i \ge 1/2.$ On the other hand, we have $X_i < 1/2+X_{i-1}$ and $X_i \le 1.$ Therefore, we have the bound 
$$X_{i-1}=\max \left(X_{i-1},\frac{1}{2} \right) \le X_i < \min \left(\frac{1}{2}+X_{i-1},1 \right)=1.$$ Furthermore, we have $X_{k+1} \ge 1/2$ and $X_{k+1} \ge X_k.$ On the other hand, we have $X_{k+1} \le \frac{1}{2} +X_k$ and $X_{k+1} \le 1.$ Therefore, we have the bound
$$\frac{1}{2}=\max \left(X_{k},\frac{1}{2} \right) \le X_{k+1} < \min \left(\frac{1}{2}+X_{k},1 \right)=\frac{1}{2}+X_{k}.$$ Moreover, for $1 < i \le k,$ we have $X_i \ge X_{i-1}.$ On the other hand, we have $X_{i} \le 1/2 +X_{i-1}$ and $X_{i} <1/2.$ Therefore, we have the bound
$$X_{i-1}=\max(X_{i-1},0) \le X_{i} < \min \left(\frac{1}2+X_{i-1},\frac{1}2 \right)=\frac{1}{2}.$$ 

Putting everything together, we see that
\begin{align}\label{p_n,k Multiple Integral}
p_{n,k} &= \int_{0}^{\frac{1}{2}} \int_{x_1}^{\frac{1}{2}} \ \dots \ \int_{x_{k-1}}^{\frac{1}{2}} \int_{\frac{1}{2}}^{\frac{1}{2} + x_k} \int_{x_{k+1}}^1 \ \dots \ \int_{x_{n-1}}^1 \ dx_n \ \dots \ dx_{k+2} \ dx_{k+1} \ dx_k \ \dots \ dx_2 \ dx_1.
\end{align}

We now evaluate the multiple integral on the right hand side of \eqref{p_n,k Multiple Integral}. First, we make the affine change of variables $x_i=\frac{1-u_i}{2}$ for $1 \le i \le k$ and $x_i=1-u_i$ for $k < i \le n,$ which has Jacobian Determinant 
$$\left | \frac{\partial(x_1, \ \dots \ , x_n)}{\partial(u_1, \ \dots \ , u_n)} \right| = \left(\frac{1}{2} \right)^k.$$ By the Change of Variables formula, we obtain

\begin{align}
p_{n,k} &= \int_{0}^{1} \int_{0}^{u_1} \ \dots \ \int_{0}^{u_{k-1}} \left( \int_{\frac{u_k}{2}}^{\frac{1}{2}} \int_{0}^{u_{k+1}} \ \dots \ \int_{0}^{u_{n-1}} \left(\frac{1}{2} \right)^k \ du_n \ \dots \ du_{k+2} \ du_{k+1} \right)  \ du_k \ \dots \ du_2 \ du_1 \label{p_n,k Multiple Integral COV} \\
&= \int_{0}^{1} \int_{0}^{u_1} \ \dots \ \int_{0}^{u_{k-1}} \frac{\left(\frac{1}{2} \right)^n}{(n-k)!} - \frac{\left(\frac{1}{2} \right)^n u_k^{n-k}}{(n-k)!} \ du_k \ \dots \ du_2 \ du_1 \label{p_n,k Multiple Integral Simplified 1} \\
&= \frac{\left(\frac{1}{2} \right)^n}{(n-k)! \ k!} - \frac{\left(\frac{1}{2} \right)^n}{(n-k)! \ (n-k+1) \ \dots \ (n) } \label{p_n,k Multiple Integral Simplified 2}   \\ 
&= \frac{\left(\frac{1}{2} \right)^n}{(n-k)! \ k!} - \frac{\left(\frac{1}{2} \right)^n}{n!} \nonumber \\
&= \frac{\binom{n}k  -1}{2^n \ n!} \nonumber,
\end{align}
where we obtained \eqref{p_n,k Multiple Integral Simplified 1} by evaluating the $n-k$ dimensional integral within the parentheses in \eqref{p_n,k Multiple Integral COV} using our integration identity from \textbf{Theorem \ref{Integration Identity of Monomial Function Over A Simplex}} and obtained \eqref{p_n,k Multiple Integral Simplified 2} by evaluating the integral in \eqref{p_n,k Multiple Integral Simplified 1} using \textbf{Theorem \ref{Integration Identity of Monomial Function Over A Simplex}} once again.
\end{proof}

We now can solve the Broken Stick problem assuming the ordering $X_1 \le \ \dots \ \le X_n$ with no extra conditions imposed. By the law of total probability, the probability that $X_1, \ \dots \ , X_n$ satisfy the conditions listed in \eqref{Broken Stick Conditions Assuming Ordering} is  equal to
\begin{align}
\sum_{k=1}^{n-1} p_{n,k} &= \sum_{k=1}^{n-1} \left(\frac{\binom{n}k  -1}{2^n \ n!}\right) \nonumber \\
&=\frac{1}{2^n \ n!} \sum_{k=1}^{n-1} \binom{n}k - \frac{1}{2^n \ n!} \sum_{k=1}^{n-1} \ 1  \nonumber \\
&= \frac{1}{2^n \ n!} \left(2^n - \binom{n}0 - \binom{n}n \right) - \frac{1}{2^n \ n!} (n-1)\label{sum p_n,k simplified}  \\
&= \frac{1}{2^n \ n!} \left(2^n - 2 \right) - \frac{1}{2^n \ n!} (n-1) \nonumber \\
&= \frac{2^n-(n+1)}{2^n \ n!} \nonumber,
\end{align} 
where we simplified the summation term in \eqref{sum p_n,k simplified} using our combinatorial identity from \textbf{Theorem \ref{Combinatorial Identity}}.

Finally, since there are $n!$ different orderings of the random variables $X_1, \ \dots \ , X_n,$ the solution to the original Broken Stick problem without any assumptions on ordering is
$$n! \left( \frac{2^n-(n+1)}{2^n \ n!} \right)= 1- \frac{n+1}{2^n}.$$

\bibliography{BrokenStickProblem.bib}
\bibliographystyle{unsrt}
\end{document}